\input amstex
%\input epsf
%\input preamble
%\newsymbol\boxtimes 1202

%letters

\def\b1{\text{\bf 1}}

\def\BC{{\Bbb C}}

\def\CD{{\Cal D}}

\def\CF{{\Cal F}}
\def\CG{{\Cal G}}

\def\CO{{\Cal O}}

\def\CT{{\Cal T}}
\def\CU{{\Cal U}}

\def\dpar{\partial}

\def\fD{{\frak D}}

\def\fG{{\frak G}}

\def\hD{{\widehat D}}

\def\hO{{{\widehat O}}}

\def\ta{\tilde a}
\def\tb{\tilde b}

\def\tg{\tilde g}
\def\tH{\tilde H}

\def\tphi{\tilde\phi}
\def\tpsi{\tilde\psi}

%symbols

\def\btu{\bigtriangleup}

\def\iso{\buildrel\sim\over\longrightarrow} 

\def\lra{\longrightarrow}

% space between paragraphs 
\parskip=6pt

\documentstyle{amsppt}
\document
%%%\magnification=1200
\NoBlackBoxes
%\nologo

%%%%%%%%%

\centerline{\bf Gerbes of chiral differential operators}

\bigskip
\centerline{Vassily Gorbounov,  
Fyodor Malikov, Vadim Schechtman}
\bigskip

\bigskip

\bigskip

{\bf 1.} In this note we compute the cohomological obstruction 
to the existence of certain sheaves of vertex algebras on smooth  
varieties. These sheaves have been introduced and studied in the previous   
work by A.Vaintrob and two of the authors, cf.  
[MSV] and [MS1]. Hopefully our result clarifies to some extent 
the constructions of {\it op. cit.} 

Recall that in [MSV] we discussed two kinds of sheaves on 
smooth complex algebraic (or analytic) varieties $X$. First, we  
defined the sheaf of conformal vertex superalgebras $\Omega^{ch}_X$, 
called {\it chiral de Rham algebra}. 
These sheaves are canonically defined 
for an arbitrary $X$. Second, for some varieties $X$ one can define a purely 
even counterpart of $\Omega^{ch}_X$, a sheaf of graded vertex 
algebras $\CO^{ch}_X$, called a {\it chiral structure sheaf}, cf. {\it 
op. cit.,} \S 5. For example, one can define $\CO^{ch}_X$ for curves, and 
for flag spaces $G/B$.  
For an arbitrary $X$, there arises certain 
cohomological obstruction to the existence of $\CO^{ch}_X$. The 
infinitesimal incarnation of this obstruction is calculated in 
{\it op. cit.}, \S 5, A.  

In the present 
note we compute explicitly this obstruction, which turns out 
to coincide with $ch_2(X)/2$,   
cf. Theorem 3 below. From a different viewpoint, this 
theorem may be regarded as a geometric interpretation of the second 
component of Chern character. In nos. 7, 8  
we give a generalization of Theorem 3, cf. Theorem 8.1. 

In no. 9 we compute the "conformal anomaly": the obstruction 
to the existence of a globally defined Virasoro field $L(z)$; 
it is given by the first Chern class $c_1(X)$, cf. Theorem 9.1.  

In particular, Theorem 3, provides a geometric criterion for a   
manifold to admit a $BU\langle 6\rangle$-structure: those are precisely 
the manifolds which admit the above mentioned sheaf $\CO^{ch}_X$ 
and for which the conformal anomaly vanishes. If such a manifold 
is Calabi-Yau (i.e. has the trivial canonical bundle) then 
$\CO^{ch}_X$ is a sheaf of {\it conformal} vertex algebras, 
cf. Corollary 9.3.       

A remark about terminology and notation. 
As was noted in [MSV], (cf. the last paragraph of \S 5), 
the sheaf $\CO^{ch}_X$ 
may be regarded as a chiral counterpart of the sheaf $\CD_X$ of 
differential operators on $X$. By that reason we prefer in this note to 
change its 
notation to $\CD^{ch}_X$ and call it {\it a sheaf of chiral differential  
operators} on $X$. From this point of view, the sheaf $\Omega^{ch}_X$ 
is a chiral counterpart of the sheaf $\CD_{\Omega_X}$ of differential 
operators on the de Rham algebra $\Omega_X$ of differential forms, 
so it would probably deserve the notation $\CD^{ch}_{\Omega_X}$, cf. 
[MS2].  

For a graded vertex algebra $V$, $V_i$ will denote the component 
of conformal weight $i$.  
In this note we will live  in the analytical category. 
 
We are deeply grateful to A.Beilinson for very 
enlightening discussions. Our special gratitude goes to him 
for an important correction in the previous version 
of the note.  

The work was done while the first and the third author visited 
the University of Glasgow and Max-Planck-Institut f\"ur Matematik; we thank 
these institutions for the hospitality.   

{\bf 2.} Fix an integer $N\geq 1$. Let $O=\BC[b^1,...,b^N]$ be the polynomial 
algebra. Following [MSV], define the vertex algebra $D$ as the free bosonic 
vertex algebra generated by the fields $b^i(z)$ of conformal weight $0$ and 
$a^i(z)$ of conformal weight $1$ ($i=1,\ldots,N$), with the OPE 
$$
a^i(z)b^j(w) \sim \frac{\delta_{ij}}{z-w},
\eqno(2.1)
$$
the other products being trivial. The algebra $D$ is called {\it the 
algebra of (algebraic) chiral differential operators on $\BC^N$}. 

The conformal weight zero component $D_0$ 
is equal to $\BC[b_0^1,\ldots,b_0^N]$, and will be identified with $O$. 
The component $D_1$ is a free $O$-module with the base $b^i_{-1}, a^i_{-1}$. 
The submodule $\Omega:=\oplus\ O\cdot b^i_{-1}\subset V_1$ will be identified 
with the $O$-module of $1$-differentials $\Omega^1_O=\oplus\ O\cdot db^i$, by 
identifying $b^i_{-1}$ with $db^i$. The submodule $T:=\oplus\ O\cdot a^i_{-1}$ 
will be identified with the $O$-module of vector fields 
$\oplus\ O\cdot\dpar_i$ where 
$\dpar_i:=\dpar/\dpar b^i$, by identifying $a_{-1}^i$ with $\dpar_i$.   

Let $X$ be an $N$-dimensional smooth complex analytic variety. 
Let $U\subset X$ be an 
open subset given together with a coordinate system, i.e. an 
open embedding $b=(b^1,\ldots,b^N):\ U\lra \BC^N$. Let $\CO(U)$ denote the 
algebra of holomorphic functions on $U$.  
The map $b$ induces an algebra homomorphism $O\lra \CO(U)$. We set 
$D(U;b):=\CO(U)\otimes_O D$. According to [MSV], 6.9, $D(U;b)$ admits 
a canonical structure of a conformal vertex algebra. It is called the 
algebra of 
(holomorphic) chiral differential operators on $U$, 
associated with the coordinate system $b$. 

If $V\subset U$ is an open subset, the map $b$ induces a coordinate system 
on $V$, and we get an algebra $D(V;b)$. When $V$ varies, we get a 
sheaf of conformal vertex algebras $\CD_{(U;b)}$ on $U$ --- the sheaf 
of holomorphic chiral differential operators on $U$, associated with 
the coordinate system $b$.   

The component $\CD_{(U;b)0}$ is identified 
with the sheaf of holomorphic functions $\CO_U$, 
and the component $\CD_{(U;b)1}$ is identified with 
the direct sum of the sheaf $\Omega^1_U$ of holomorphic $1$-differentials  
and the sheaf $\CT_U$ of holomorphic vector fields. 

Let us define a two-step filtration on $\CD_{(U;b)1}$ by 
$$
0=F^2\CD_{(U;b)1}\subset \Omega^1_U=F^1\CD_{(U;b)1}\subset 
\CD_{(U;b)1}=F^0\CD_{(U;b)1}
\eqno{(2.2)}
$$
Obviously, $F^0/F^1=\CT_U$. The associated graded space 
$$
gr^._F\CD_{(U;b)1}=\Omega^1_U\oplus \CT_U
\eqno{(2.3)}
$$
does not depend on the choice of a coordinate system $b$, cf. 
[MSV], 6.2.  

Let us define a groupoid $\fD^{ch}(U)$ as follows. The objects of $\fD^{ch}(U)$ 
are coordinate systems $b$ on $U$. Given two coordinate systems 
$b,\ 'b$, the set of morphisms $Hom_{\fD^{ch}(U)}(b,\ 'b)$ consists of 
isomorphisms of sheaves of graded vertex algebras 
$\CD_{(U;b)}\iso \CD_{(U;'b)}$ 
which induce the identity on the conformal weight zero components, 
and on the associated graded (with respect to the filtration (2.2)) 
of the conformal weight one components. (Note that any morphism of graded  
vertex algebras as above necessarily respects the filtration (2.2)). 
This groupoid will be called the {\it groupoid of (holomorphic) chiral 
differential 
operators} on $U$. When $U$ varies, we get a stack ({\it champ}) 
of groupoids, i.e. a {\it gerbe} $\fD^{ch}_X$ on $X$, called 
the {\it gerbe of (holomorphic) chiral differential operators} on $X$.

Concerning the J.Giraud's language of gerbes, see for example 
Section 5 of [D]. 
Now we can formulate our main result.
 
%\newpage

{\bf 3. Theorem.} {\it The} lien {\it of the gerbe 
$\fD^{ch}_X$ is canonically 
isomorphic to the sheaf of closed holomorphic $2$-differentials 
$\Omega^{2, cl}_X$. 

The equivalence class $c(\fD^{ch}_X)\in 
H^2(X;\Omega^{2, cl}_X)$ is equal to 
$$
c(\fD^{ch}_X)=2 ch_2(\CT_X)
\eqno{(3.1)}
$$} 

Here $\CT_X$ denotes the tangent bundle; for a vector bundle $E$ we set
$$
ch_2(E)=\frac{1}{2}c_1^2(E)-c_2(E)
\eqno{(3.2)}
$$ 
where 
$c_i(E)\in H^i(X;\Omega^{i, cl}_X)$ are the Hodge Chern classes.
 
Recall that the set of equivalence classes of gerbes on $X$ bounded 
by a given abelian sheaf $\CF$ is identified with the group $H^2(X;\CF)$. 
Namely, if $\fG$ is such a gerbe, choose an open covering $\CU=\{ U_i\}$ 
of $X$  
such that for each $U_i$ the groupoid $\fG(U)$ is nonempty and connected. 
Choose objects $\CG_i\in\fG(U_i)$. On the pairwise intersections there exist 
isomorphisms
$\phi_{ij}:\ \CG_i|_{U_{ij}}\iso\ \CG_j|_{U_{ij}}$. The "discrepancies" 
$c_{ijk}:=\phi_{ij}\phi_{jk}\phi_{ik}^{-1}$ (on triple intersections) 
define a Cech $2$-cocycle 
in $Z^2(\CU;\CF)$ which represents the class $c(\fG)$. 
 
This theorem means that a "sheaf of chiral differential operators"  
is not a uniquely defined object associated with $X$, but an object 
of a canonically defined {\it groupoid} $\Gamma(X;\fD^{ch}_X)$. 
This groupoid is empty if the class (3.1) is non-trivial. 

If (3.1) is trivial then this groupoid is equivalent 
to the 
groupoid of $\Omega^{2, cl}_X$-torseurs on $X$. Thus, in this case the set 
of isomorphism classes of $\Gamma(X;\fD^{ch}_X)$ is  isomorphic to 
$H^1(X;\Omega^{2, cl}_X)$, and the automorphism group of an object 
is isomorphic to $H^0(X;\Omega^{2, cl}_X)$. 

 From this one sees immediately that for curves the obstruction vanishes, and  
the above groupoid 
is trivial: one can define canonically "the" sheaf of chiral differential 
operators $\CD_X^{ch}$. For the explicit construction, see (6.2) below.   

The same is true for the flag spaces $G/B$. 
This explains the constructions of [MSV]. As a negative example, note that 
for the projective spaces of dimension greater than $1$, the class 
(3.1) is non-zero. 

{\bf 4.} Let us explain a local statement which implies Theorem 3. 
Let $\hO\supset O$ be the algebra of the Taylor power series in variables 
$b^i$ convergent in some neighbourhood of the origin. Let $G$ be the group 
whose elements are $N$-tuples $g=(g^1(b),\ldots,g^N(b))$ where $g^i\in \hO$ and   
the Jacobian matrix $\dpar g:=(\dpar_ig^j)$ is non-degenerate. 
The composition is defined by 
$$
(g_1 g_2)^i(b)=g_2^i(g_1^1(b),\ldots,g_1^N(b))
\eqno{(4.1)}
$$
The group $G$ acts in the obvious way from the right on the algebra $\hO$. 

We have 
$$
\dpar(g_1 g_2)=\dpar g_1 \dpar g_{2 g_1}
\eqno{(4.2)}
$$
 Here we use the notation 
$$              
h_g(b):= h(g(b))
\eqno{(4.3)}
$$
We will denote by $\Omega^1_{\hO}$ (resp., $\CT_{\hO}$) the module 
of germs of holomorphic $1$-differentials (resp., holomorphic vector 
fields). It is a free $\hO$-module with the base $\{ db^i\}$ 
(resp., $\{\dpar_i\}$). 

Let us consider the vertex algebra $\hD:=\hO\otimes_O D$. We have 
$\hD_0=O$, and $\hD_1$ has a canonical two step filtration defined as 
in (2.2), with $F^1\hD_1=\Omega^1_{\hO},\ \hD_1/F^1\hD_1=\CT_{\hO}$. 

Let us call an automorphism of $\hD$ natural if it induces an automorphism 
$g\in G$ on $\hD_0$, and the induced by $g$ automorphism on $gr^._F\hD_1$. 

We want to describe the group $G^\sim$ of all natural automorphisms of 
$\hD$. In coordinates, a natural automorphism is given by a pair 
$(g,h)$ where $g\in G$ and $h\in M_N(\hO)$ is an $N\times N$-matrix 
over $\hO$. By definition, such a pair acts on the generating fields as 
follows
$$
(b^i\circ (g,h))(z)=g^i(b)(z)
\eqno{(4.4a)}
$$
$$
(a^i\circ (g,h))(z)=a^j\phi^{ji}(g)(z)+b^k(z)'h^{ki}(z)
\eqno{(4.4b)}
$$
(this will be the right action). Here the matrix $\phi(g)=(\phi^{ij}(g))$ 
is defined by 
$$
\phi(g)=(\dpar g)^{-1 t}
\eqno{(4.5)}
$$
where $(a^t)^{ij}=a^{ji}$. 
  
The straightforward computation of OPE shows that a pair $(g,h)$ defines  
an automorphism of the vertex algebra $\hD$ if and only if the matrix 
$h$ satisfies the equations 
$$
h^t\phi(g)+\phi^t(g) h=\psi(g)
\eqno{(4.6a)}
$$
and
$$
\phi^{pi}\dpar_p h^{sj} - \phi^{pj}\dpar_p h^{si} +
$$
$$ 
+ \dpar_s \phi^{pi} h^{pj} + \dpar_s h^{pi} \phi^{pj} - 
\dpar_s\dpar_r\phi^{pi}\dpar_p\phi^{rj}=0
\eqno{(4.6b)}
$$
where $\psi(g)=(\psi^{ij}(g))$ is a symmetric matrix 
$$
\psi^{ij}(g)=\dpar_r\phi^{pi}(g)\dpar_p\phi^{rj}(g)
\eqno{(4.7)}
$$
(we imply everywhere the summation over the repeating indices). 

The computation similar to that in the proof of Theorem 3.8 from  
[MSV], shows that the composition on $G^\sim$ 
is given by the formula
$$
(g_1,h_1)(g_2,h_1)=(g_1g_2, h_1\phi(g_2)_{g_1}+\dpar g_1 h_{2 g_1} + 
\alpha(g_1,g_2))
\eqno{(4.8)}
$$
where $\alpha(g_1,g_2)=(\alpha^{ij}(g_1,g_2))$ is given by 
$$
\alpha^{ij}(g_1,g_2)=\dpar_i\phi^{pq}(g_1)\dpar_p[\phi^{qj}(g_2)_{g_1}]
\eqno{(4.9)}
$$
This describes the group $G^\sim$. The unit is $(1,0)$. Let us introduce  
a bigger group $G'$ consisting of all pairs $(g,h)$ as above, with 
the multiplication (4.8), but with $h$ not necessarily satisfying 
the conditions (4.6 a, b). One has to check 
the associativity of multiplication, which is straightforward. 
Obviously, $G^\sim\subset G'$.  

It is convenient to introduce in $G'$ the new coordinates by setting 
$$
(g,h)^\sim := (g,0)(1,h)=(g,h\dpar g^{-1 t})
\eqno{(4.9)}
$$
In new coordinates the multiplication in $G'$ looks as follows 
$$
(g_1,h_1)^\sim (g_2,h_2)^\sim = 
(g_1g_2, h_1 + \dpar g_1 h_{2 g_1} \dpar g_1^t +\beta (g_1,g_2))^\sim
\eqno{(4.10)}
$$
where 
$$
\beta(g_1,g_2)=\alpha(g_1,g_2)\dpar(g_1g_2)^t
\eqno{(4.11)}
$$
In particular, 
$$
(1,h)(g,0)=(h,\dpar g\cdot h_g\cdot \dpar g^t)
\eqno{(4.12)}
$$
We have an obvious epimorphism $G'\lra G$, $(g,h)\mapsto g$, and the formula 
(4.12) shows that its kernel, as a $G$-module, may be identified with the 
$G$-module of quadratic differentials $\Omega^{1\otimes 2}_{\hO}:=
\Omega^1_{\hO}\otimes_{\hO} \Omega^1_{\hO}$, i.e. we have the group extension 
$$
0\lra \Omega^{1\otimes 2}_{\hO}\lra G'\lra G\lra 1
\eqno{(4.13)}
$$
Equation (4.10) means that if we chose a section of the epimorphism 
$G'\lra G$ as $g\mapsto (g,0)$ then the extension (4.13) corresponds 
to the $2$-cocycle $\beta(g_1,g_2)\in Z^2(G;\Omega^2_{\hO})$. 
The equation of a cocycle (equivalent to the associativity of multiplication 
in $G'$) reads as 
$$
\beta(g_1,g_2)-\beta(g_1,g_2g_3)+\beta(g_1g_2,g_3)-g_1^*\beta(g_2,g_3)=0
\eqno{(4.14)}
$$ 
A little computation using the formula $\dpar_i(A^{-1})= - A^{-1}\dpar_i A 
\cdot A^{-1}$ gives the following expression for $\beta$: 
$$
\beta^{ij}(g_1,g_2)=\dpar (g_1^{-1})^{qa} \dpar_i\dpar_a g_1^b 
\dpar ^j [ (\dpar_b g_2^r)_{g_1} ] (\dpar g_2^{-1})_{g_1}^{ rq}
\eqno{(4.15)}
$$
which may be rewritten in the following compact form: 
$$
\beta(g_1,g_2)=tr \{ \dpar g_1^{-1} d\dpar g_1 \otimes g_1^* 
(d \dpar g_2\cdot \dpar g_2^{-1}) \} 
\eqno{(4.16)}
$$

{\bf 5.} Now let us take care of the subgroup 
$G''\subset G'$ 
given inside $G'$ by the equation (4.6a). Since the matrix $\psi(g)$ is 
symmetric, and the matrix $h$ is arbitrary, the above equation has 
$N(N-1)/2$-dimensional (over $\hO$) space of solutions. 

For example, for $N=1$, (4.6) has a unique solution, 
and $G'' =G$, namely
$$
h(g)=g^{\prime\prime 2}/2g^{\prime 3}
\eqno{(5.1)}
$$ 
 
%every automorphism $g\in G$ lifts uniquely to a natural automorphism of 
%$\hD$ given by 
%$$
%\tilde{g}=(g,g^{\prime\prime 2}/2g') 
%\eqno{(5.1)}
%$$
%This corresponds to the fact that for curves the sheaf $\CD^{ch}_X$ 
%is uniquely defined. Note that the formula [MSV], (5.23b) is a particular 
%case of (5.1) (for $g(b)=b^{-1}$). 

Now let us treat the case of an arbitrary $N$. First of all, $\psi(1)=0$ 
and $\phi(1)=1$, and the equation (4.6) means simply that the matrix $h$ 
is skew symmetric. This means that $G''$ is included in 
a group extension 
$$
0\lra \Omega^2_{\hO}\lra G'' \lra G\lra 1
\eqno{(5.2)}
$$ 
It is easy to write down a particular solution of (4.6a),  
namely
$$
h(g)=\phi^t(g)^{-1}\psi(g)/2
\eqno{(5.3)}
$$
This gives a section of (5.2), given by 
$$
g\mapsto (g,h(g))=: (g,s(g))^\sim
\eqno{(5.4)}
$$
If $c(g_1,g_2)\in Z^2(G;\Omega^2_{\hO})$ is the cocycle corresponding to this 
choice of the section then we have 
$$
c(g_1,g_2)=\beta(g_1,g_2)+s(g_1)-s(g_1g_2)+\dpar g_1 s(g_2)_{g_1} \dpar g_1^t
\eqno{(5.5)}
$$
 In other words, $c$ is cohomologous to $\beta$. But as far as we know 
 that the cohomology class of $\beta$ comes from a cohomology class 
 from $H^2(G;\Omega^2_{\hO})$, this last class is unique and may be given 
 simply by the skew symmetrization of $\beta$, i.e. 
 by 
 $$
 \gamma(g_1,g_2)= tr\{ \dpar g_1^{-1} d \dpar g_1\wedge 
 g_1^* (d \dpar g_2 \cdot \dpar g_2^{-1})\}=
 $$
 $$
 =tr\{d \dpar g_1\cdot \dpar g_1^{-1}\wedge \dpar g_1 \cdot g_1^*[ 
 d \dpar g_2\cdot \dpar g_2^{-1}]\cdot \dpar g_1^{-1}\} 
 \eqno{(5.6)}
 $$
 Comparing the expression (5.6) with the formula (4.5) from [TT] we see that 
 the cocycle $\gamma$ is nothing but the universal cocycle of $G$ 
 representing the characteristic class corresponding to the invariant 
 polynomial $tr(AB)$ (on pairs of $N\times N$ matrices), i.e. 
 it represents the component of the Chern character $c_1^2-2c_2$. 
 
 Thus, we arrive at 
 
 {\bf 5.1. Proposition.} {\it The group $G''$ is 
 included in the group extension 
 $$
 0\lra \Omega^2_{\hO}\lra G''\lra G\lra 1
 \eqno{(5.7)}
 $$
 The cohomology class of this extension is equal to 
 $2ch_2:=c_1^2-2c_2\in H^2(G;\Omega^2_{\hO})$ where $c_i\in H^i(G;\Omega^i_{\hO})$ 
 are the universal Chern classes.} 

{\bf 6.} Finally, group $G^\sim$ is given inside $G''$ by 
the equation (4.6b). For $g=1$ this equation simply means 
that the form $h^{ij}db^i\wedge db^j$ is closed 
(we are grateful to A.Beilinson for this remark). 

Therefore, the kernel of the obvious map 
$G^\sim\lra G$ is isomorphic to the $G$-module 
$\Omega_{\hO}^{2, cl}$ of closed $2$-forms.

{\bf 6.1. Lemma.} {\it The function $h(g)$ given by the formula (5.3) 
satisfies the equation (4.6b).} 

In fact, if we symmetrize the equation (4.6b) with respect to $i, j$, 
we get $\dpar_s$(4.6a). Therefore, it is enough to check the 
skew symmetrized part of (4.6b). This is done by a direct computation, using 
(4.5). 

This lemma gives a section of the canonical map $G^\sim\lra G$; in particular 
the last map is surjective. 
 
The above considerations, together with Proposition 5.1, imply   

{\bf 6.2. Theorem.} {\it The group $G^\sim$ of all natural 
automorphisms of the vertex algebra $\hD$ is included in 
the group extension
$$
0\lra \Omega^{2,cl}\lra G^\sim \lra G\lra 1
\eqno{(6.1)}
$$
The cohomology class of this extension is equal to 
$2ch_2:=c_1^2-2c_2\in H^2(G;\Omega^{2, cl})$ where 
$c_i\in H^i(G;\Omega^{i, cl})$ are the universal 
Chern classes.} 

Note that the standard cocycles representing the universal 
Chern classes in $H^i(G;\Omega^i_{\hO})$ in fact 
take values in 
$\Omega^{i, cl}_{\hO}$. 

Let us sketch another argument proving surjectivity of   
the map $G^\sim\lra G$.  
Note that in [MSV], \S 5 one has 
lifted an arbitrary element of the Lie algebra 
$W_N:=Lie(G)$ to a derivation $\pi(\tau)$ of $\hD$. 
Since every 
$g\in G$ is the exponent of some vector field $\tau$, 
we can take $\tilde{g}=\exp (\pi(\tau))$ as a lifting 
of $g$. 

Consider the case $N=1$. In this case $G^\sim=G''=G$. Thus,  
every automorphism $g\in G$ lifts uniquely to a natural automorphism of 
$\hD$ given by 
$$
\tilde{g}=(g,g^{\prime\prime 2}/2g^{\prime 3}) 
\eqno{(6.2)}
$$
cf. (5.1). 
This corresponds to the fact that for curves the sheaf $\CD^{ch}_X$ 
is uniquely defined. Note that the formula [MSV], (5.23b) is a particular 
case of (6.2) (for $g(b)=b^{-1}$).

Our main Theorem 3 is an immediate consequence of 
Theorem 6.1. 
 
 {\bf 7.} Let $E$ be a vector bundle over $X$, $\Lambda E$ its 
 exterior algebra. Proceeding in a similar way as above, we can 
 produce a gerbe $\fD^{ch}_{\Lambda E}$ of chiral differential 
 operators on $\Lambda E$, also bounded by the lien $\Omega^{2, cl}_X$. 

 Let us explain how to do this. Let us describe the local model. 
 In the notations of no. 2, fix an integer $M\geq 0$.  
 Let us define a vertex superalgebra 
 $D_\Lambda\supset D$ which is generated by the free bosonic fields 
 $b^i(z), a^i(z)$ as in {\it loc. cit.} and free odd fields 
 $\phi^j(z)$ of conformal weight $0$ and $\psi^j(z)$ of conformal 
 weight $1\ (j=1,\ldots, M)$, with the OPE (2.1) and 
 $$
 \psi^i(z)\phi^j(w)\sim \frac{\delta_{ij}}{z-w}
 \eqno{(7.1)}
 $$ 
 the other products being trivial. 
Set $\hD_\Lambda:=\hO\otimes_O D_\Lambda$. Consider the group $H$ which 
is the semidirect product of the group of coordinate changes $G$ from no. 4, 
the group of $M\times M$ matrices $Mat_M(\hO)$, with the obvious action 
of $G$ on $Mat_M(\hO)$. The group $H$ will be the analogue of the group 
$G$ in our situation. 

Let us call an automorphism of the superalgebra $\hD_\Lambda$ 
{\it admissible} if it has the form
$$
\tb^i(z)=g^i(b)(z)
\eqno{(7.2a)}
$$
$$
\ta^i(z)=a^jc^{ji}(g)(z) -\psi^p\phi^q A^{qs} c^{ki}\dpar_k A^{-1 sp}(z) +
 b^j(z)' h^{ij}(z)
\eqno{(7.2b)}
$$
$$
\tphi^i(z)=\phi^j A^{ji}(z)
\eqno{(7.2c)}
$$
$$
\tpsi^i(z)=\psi^j A^{-1 ij}(z)
\eqno{(7.2d)}
$$
for some $(g, A)\in H$ and $h\in Mat_N(\hO)$. Here $c(g):=(\dpar g)^{-1 t}$ 
(we have changed the notation from (4.5), since we already have $\phi$'s).  

We want to compute the group $\tH$ of all admissible automorphisms of 
$\hD_\Lambda$. The composition will be 
$$
(g_1,A_1,h_1)(g_2,A_2,h_2)=
$$
$$
=(g_1g_2,A_1A_{2g_1},h_1c(g_2)_{g_1}+\dpar g_1h_{2g_1}+\alpha(g_1,g_2) -
\gamma(g_1,g_2))
\eqno{(7.3)}
$$
where $\alpha(g_1,g_2)$ is given by (4.9), and 
$$
\gamma(g_1,A_1;g_2,A_2)^{ij}=
\dpar_i A_1^{-1 pr} A_1^{rq} A_{2 g_1}^{qs} c(g_2)_{g_1}^{k j}\dpar _k 
A_{2 g_1}^{-1 sp}
\eqno{(7.4)}
$$
Set
$$
\delta(g_1,A_1;,g_2,A_2)=\gamma(g_1,A_1,g_2,A_2)\dpar(g_1g_2)^t
\eqno{(7.5)}
$$
cf. (4.11). A little computation gives  

{\bf 7.1. Lemma.} {\it We have
$$
\delta(g_1,A_1,g_2,A_2)=tr\{ dA_1\cdot A_1^{-1}\otimes 
A_1 g_1^* (dA_2\cdot A_2^{-1}) A_1^{-1}\}
\eqno{(7.6)}
$$}
Cf. (4.16). 

Let $H'$ be the group of all triples $(g,A,h)$ as above, with the composition 
given by (7.3). We arrive at 

{\bf 7.2. Proposition.} {\it The group $H'$ is included in the extension
$$
0\lra \Omega_{\hO}^{ 1\otimes 2}\lra H'\lra H\lra 1
\eqno{(7.7)}
$$
given by the two-cocycle 
$$
\epsilon(g_1,A_1;g_2,A_2)=
\beta(g_1,g_2)-\gamma(g_1,A_1;g_2,A_2)
\eqno{(7.8)}
$$
where $\beta$ (resp. $\gamma$) is given by (4.16) (resp. by (7.6)).} 

Now, proceeding in the same way as in nos. 5, 6, we arrive at 

{\bf 7.3. Proposition}. {\it The group $\tH$ of all admissible 
automorphisms of the vertex algebra $\hD_\Lambda$ is included into 
the extension
$$
0\lra \Omega^{2, cl}_{\hO}\lra\tH\lra H\lra 1
\eqno{(7.9)}
$$
whose cohomology class is given by the skew symmetrization of the cocycle 
(7.8), i.e. equal to 
$$
c_1^2-2c_2 - (c_1^{\prime 2}-2c'_2)\in H^2(H,\Omega^{2, cl}_{\hO})
\eqno{(7.10)}
$$
where $c_i, c'_i\in H^i(H;\Omega^{i, cl}_{\hO})$ are the universal 
Chern classes.}

{\bf 8.} Returning to the global situation of the beginning of the 
preceeding subsection (we assume that $dim(X)=N,\ rk(E)=M$), we 
define for a sufficiently small $U\subset X$ the groupoid 
$\fD^{ch}_{\Lambda E}(U)$ whose objects are pairs 

(a coordinate system on $U$, a trivialization of $E$ over $U$), 

and morphisms are defined as 
in no. 2, using the admissible isomorphisms of the corresponding 
vertex superalgebras. This way we get the gerbe $\fD^{ch}_{\Lambda E}$ over $X$. 

Proposition 7.3 implies 

{\bf 8.1. Theorem.} {\it The gerbe $\fD_{\Lambda E}^{ch}$ is bounded by the 
lien $\Omega^{2, cl}_X$, and its equivalence class is equal to 
$$
c(\fD_{\Lambda E}^{ch})=2ch_2(\CT_X)-2ch_2(E)\in H^2(X,\Omega^{2,cl}_X)
\eqno{(8.1)}
$$}

Note that the class (8.1) does not change when $E$ is replaced by the dual 
bundle. 

We see that for $E=\Omega^1_X$, or $\CT_X$ we arrive 
at the neutral gerbe. In fact, the groupoid of global sections 
$\fD^{ch}_{\Lambda E}(X)$ in both cases admits a canonical object. 
For $E=\Omega^1_X$ it is the "chiral de Rham complex" from 
[MSV], and for $E=\CT_X$ it is a sort of a "mirror partner". 

{\bf 9. conformal anomaly.} Let us return to the situation of 
nos. 4 --- 6. Consider an automorphism of $\hD$ given by an 
element $\tg=(g,h(g))$ where $h(g)$ is given by (5.3).    
Recall that $\hD$ admits a canonical structure of a {\it conformal} 
vertex algebra, with the Virasoro field 
$$
L(z)=a^i(z)b^i(z)'
\eqno{(9.1)}
$$
{\bf 9.1. Theorem.} {\it Under the automorphism $\tg$, the Virasoro 
field $L(z)$ is transformed as follows
$$
L(z)^\sim=L(z)-\frac{1}{2}(Tr\ \log\ \dpar g(b(z)))^{\prime\prime}
\eqno{(9.2)}
$$} 

Here the tilde denotes the result of transformation under $\tg$. 

Compare this formula with [MSV], Theorem 4.2.  

{\bf Proof.} Let us consider the one-dimensional case for simplicity. 
We have, in the obvious notations, 
$$
L^\sim=a(z)^\sim_{-1}b(z)^\sim_{-1}\b1=
g(b(z))_{-1}\bigl[a_{-1}g'_b(b_0)^{-1}+b_{-1}
\frac{g^{\prime\prime}_{bb}(b_0)^2}{2g'_b(b_0)^3}\bigr]\b1
\eqno{(9.3)}
$$
We have the Taylor formula
$$
g(b(z))=g(b_0)+g'_b(b_0)\Delta b(z)+\frac{1}{2}g''_{bb}(b_0)
\Delta b(z)^2+\frac{1}{6}g'''_{bbb}(b_0)\Delta b(z)^3+\ldots
\eqno{(9.4)}
$$
where 
$$
\Delta b(z)=b(z)-b_0=b_{-1}z+b_1z^{-1}+b_{-2}z+b_2z^{-2}+\ldots
$$
The component $g(b(z))_{-1}$ is the coefficient at $z$ of the expression 
(9.4). This is of course the infinite sum, but we are interested 
only in the terms which give the nonzero contribution after 
plugging into (9.3). These are the terms containing (except for $b_0$) 
the monomials  
$b_1$, $b_{-1}^2b_1$ and $b_{-2}b_1$ (the last two would interact with 
$a_{-1}$ in (9.3)). Let us write down these terms:
$$
g(b(z))_{-1}=g'_b(b_0)b_{-1} +
g''_{bb}(b_0)b_{-2}b_1+
\frac{1}{2}g'''_{bbb}(b_0)b_{-1}^2b_1+\ldots
$$
Substituting this into (9.3), we get
$$
L^\sim=L-\frac{1}{2}\bigl[(2g''_{bb}(b_0)g'_b(b_0)b_{-2}+
g'''_{bbb}(b_0)g'_b(b_0)b_{-1}^2-
g''_{bb}(b_0)^2b_{-1}^2)g'_b(b_0)^{-2}\bigr]\b1
\eqno{(9.5)}
$$
Let $L_{-1}:\ \hD\lra \hD$ denote the canonical derivation of 
our vertex algebra. It remains to notice that the anomalious term 
in (9.5) is equal to 
$$
-\frac{1}{2}L_{-1}\bigl\{g''_{bb}(b_0)g'_b(b_0)^{-1}b_{-1}\b1\bigr\}
$$ 
This follows from the Leibniz rule and the formulas 
$L_{-1}(b_0\b1)=b_{-1}\b1;\ L_{-1}(b_{-1}\b1)=2b_{-2}\b1$. 
Therefore, the corresponding field is 
$$
-\frac{1}{2}\bigl\{g''_{bb}(b(z))b(z)'g'_b(b(z))^{-1}\bigr\}'= 
-\frac{1}{2}\bigl\{\log g'_b(b(z))\bigr\}''
$$
This proves the theorem in the one-dimensional case. 

The higher dimensional case is checked in a similar manner. 
We leave it to the reader. $\btu$

{\bf 9.2.} Note that (9.2) implies that the Fourier components 
$L_0,\ L_{-1}$ are always conserved as it should be: our 
sheaves $\CD^{ch}_X$ are always the sheaves of graded vertex algebras. 

{\bf 9.3. Corollary.} {\it Let $X$ be smooth variety with the trivial 
canonical bundle and with $c_2(\CT_X)=0$. Then each sheaf 
$\CD^{ch}_X\in\fD^{ch}_X(X)$ admits a structure of a sheaf of conformal 
vertex algebras.} 

Indeed, since the canonical bundle is trivial, we can choose a 
coordinate atlas of $X$ 
such that the Jacobians of the transition functions are constants. 
For such an atlas, the anomalious terms in (9.2) disappear. 
If $c_2(\CT_X)=0$ then we glue our sheaves $\CD^{ch}_X$ using the 
gluing automorphisms (4.4) were the matrices $h$ can be written as 
$h=h_1+h_2$ where $h_1$ is given by (5.3) and $h_2$ is a skew 
symmetric matrix. Note that a correction by a skew symmetric matrix 
does not change the Virasoro field: we have  
$$
b^i(z)'\bigl[ a^i(z)+b^k(z)'h^{ki}(z)\bigr]=
b^i(z)'a^i(z)+b^i(z)'b^k(z)'h^{ki}(z)=
b^i(z)'a^i(z)
$$
if $(h^{ki})$ is skew symmetric. Therefore, the Virasoro field 
is conserved in our sheaf.   

{\bf 9.4.} Similarly to nos. 7, 8, we may add the exterior bundle 
to the picture. The corresponding conformal anomaly will be equal 
$$
\frac{1}{2}(c_1(\CT_X)-c_1(E))
$$
We leave the details to the reader.

\bigskip
\bigskip 

\centerline{\bf References}

\bigskip
\bigskip 

%%%%%%
[D] 
[ 

%[K] V.~Kac, Vertex algebras for beginners, University Lecture Series, 
%{\bf 10}, American Mathematical Society, Providence, RI, 1997.

[D] P.~Deligne, Le symbole mod\'er\'e, {\it Publ. Math. IHES}, {\bf 73} (1991), 
147 -181.  

[MSV] F.~Malikov, V.~Schechtman, A.~Vaintrob, Chiral de Rham complex,\ 
{\it Comm. Math. Phys.} (1999), to appear; math.AG/9803041. 
 
[MS1] F.~Malikov, V.~Schechtman, Chiral de Rham complex. II, 
{\it D.B.~Fuchs' 60-th Anniversary volume} (1999), to appear; 
math.AG/9901065. 

[MS2] F.~Malikov, V.~Schechtman, Chiral Poincar\'e duality, math.AG/9905008. 

[TT] D.~Toledo, Y.L.~Tong, A parametrix for $\bar{\dpar}$ and Riemann-Roch 
in Cech theory, {\it Topology}, {\bf 15} (1976), 273 - 301.

\bigskip

\bigskip

V.G.: Department of Mathematics, University of Kentucky, 
Lexington, KY 40506, USA;\ vgorb\@ms.uky.edu

F.M.: Department of Mathematics, University of Southern California, 
Los Angeles, CA 90089, USA;\ fmalikov\@mathj.usc.edu   

V.S.: Department of Mathematics, University of Glasgow, 
15 University Gardens, Glasgow G12 8QW, UK;\ 
vs\@maths.gla.ac.uk,\ vadik\@mpim-bonn.mpg.de

\enddocument